\documentclass{amsart}
\usepackage{amssymb}

\usepackage{amscd}
\usepackage{thmdefs}

\input{tcilatex}
\theoremstyle{definition}
\theoremstyle{remark}
\numberwithin{equation}{section}

\input tcilatex

\begin{document}
\title[Boundary Quotient Graphs and the Graph Index]{Boundary Quotient Graphs and the Graph Index}
\author{Ilwoo Cho}
\address{Saint Ambrose University, Dep. of Math, Davenport, U. S. A.}
\email{choillwoo@sau.edu}
\date{}
\dedicatory{{\small (Communicated with Prof. Thomas Anderson)}}
\thanks{I really appreciate to all supports from the department of Mathematics,
Saint Ambrose University.}
\maketitle

\begin{abstract}
In this paper, we will consider the boundary quotient graphs. Let $G$ be a
finite directed graph with its vertex set $V(G)$ and its edge set $E(G).$
The boundary $\partial $ of $G$ is a subset of the vertex set $V(G).$ For
the given boundary $\partial $ $\subseteq $ $V(G),$ we give an boundary
quotient : if $v_{1},$ $v_{2}$ $\in $ $\partial ,$ then $v_{1}$ $=$ $v_{2},$
for all $v_{1},$ $v_{2}$ $\in $ $\partial $ Then we can construct a new
graph $G_{\partial }$ $=$ $G$ $/$ $\partial $ called the $\partial $%
-quotient graph of $G.$ In Chapter 1, we restrict our interests to the
finite simplicial directed graphs. We will observe some properties of $G$ $/$
$\partial .$ In particular, we show that all total boundary quotient graphs
has the same type, where total boundary $\partial $ is $V(G).$ Every total
boundary quotient graph is graph-isomorphic to one-vertex-$\left|
E(G)\right| $-loop-edge graph. In fact, every total boundary quotient graph
of a finite directed graph is graph-isomorphic to the
one-vertex-multi-loop-edge graph. This result shows that boundary quotient $%
\partial $ is not an invariants on finite simplicial directed graphs.
However, we show that the ``admissible'' boundary quotient is an invariant
on finite simplicial directed graphs with mixed maximal types. In Chapter 2,
we consider arbitrary finite directed graphs and define the subgraph
boundary quotient $\partial _{H}$ of the given graph $G,$ where $H$ is a
full subgraph of $G.$ The subgraph boundary $\partial _{H}$ is defined by
the set $V(H)$ $\cup $ $E(H).$ By identifying all element in $\partial _{H}$
in $G,$ we can get the subgraph boundary quotient graph $G$ $/$ $\partial
_{H}.$ Define the subgraph index $ind_{G}(\partial _{H})$ of $G$ with
respect to $H,$ by the exponetial of $V(G$ $/$ $\partial _{H})$ $\cup $ $E(G$
$/$ $\partial _{H}).$ We will observe the properties of the subgraph
boundary index of finite directed graphs. In particular, we can get that $%
Ind_{G}(\partial _{H})$ $=$ $\frac{Ind_{G}(1)}{Ind_{H}(1)},$ where $1$ is
the trivial graph.
\end{abstract}

\strut

In this paper, we will consider boundary quotient graphs. Let $G$ be a
finite directed graph with its vertex set $V(G)$ and its edge set $E(G).$ A
boundary $\partial $ of the graph $G$ is defined by a subset of $V(G).$ For
the fixed boundary $\partial $ of $G,$ the boundary quotient, also denoted
by $\partial ,$ is defined by the following relations;

\strut

(i) \ \ All elements in $\partial $ are identified. i.e., if $v_{1}$ $\neq $ 
$v_{2}$ $\in $ $\partial ,$ then we assume that $v_{1}$ $=$ $v_{2},$ for all 
$v_{1},$ $v_{2}$ $\in $ $\partial .$ In other words, the boundary $\partial $
makes a base point $v_{\partial },$ as the identified vertices by $\partial
. $

\strut

(ii) \ If $v_{1}$ $\neq $ $v_{2}$ in $\partial $ and if there exists an edge 
$e$ connecting $v_{1}$ and $v_{2},$ then the edge $e$ is replaced by a
loop-edge concentrated on the base point $v_{\partial }.$

\strut

The loop-edges in (ii) are called the $\partial $-loop-edges at the base
point $v_{\partial }.$ Also, the edges in $E(G),$ which are not affected by
the boundary quotient $\partial $, are said to be non-$\partial .$ We define
a boundary quotient graph $G$ $/$ $\partial $ of $G$ by $\partial ,$ by the
directed graph $Q$ with its vertex set

\strut

\begin{center}
$V(Q)=\{v_{\partial }\}\cup \{v\in V(G):v\notin \partial \}$
\end{center}

\strut

and its edge set

\strut

\begin{center}
$
\begin{array}{ll}
E(Q)= & \{e\in E(G):e\text{ is non-}\partial \} \\ 
&  \\ 
& \text{ }\cup \{l:l\text{ is a }\partial \text{-loop-edges at }v_{\partial
}\}.
\end{array}
$
\end{center}

\strut

We say that the boundary $\partial $ is total, if $\partial $ $=$ $V(G).$

\strut

In Chapter 1, we will consider some properties of boundary quotient graph $Q$
$=$ $G$ $/$ $\partial .$ We will show that the total boundary quotient graph 
$G$ $/$ $\partial _{t}$ of a finite simplicial directed graph $G$ is
graph-isomorphic to the one-vertex-$\left| E(G)\right| $-loop-edge graph $%
L_{\left| E(G)\right| },$ where $\partial _{t}$ is the total boundary of $G.$
This shows us that the boundary quotient is generally not an invariant on
finite simplicial directed graphs. For example, if $C_{N}$ is a one-flow
circulant graph and if $T_{N+1}$ is a finite directed tree with $(N$ $+$ $1)$%
-vertices, then clearly graphs $C_{N}$ and $T_{N+1}$ is not
graph-isomorphic. However, the total boundary quotient graphs $C_{N}$ $/$ $%
\partial _{C_{N}}$ and $T_{N+1}$ $/$ $\partial _{T_{N+1}}$ are
graph-isomorphic, because both of them are graph-isomorphic to the graph $%
L_{N},$ where $L_{N}$ is the one-vertex-$N$-loop-edge graph.

\strut

We also define so-called the admissible boundaries and the admissible
boundary quotient graphs. Let $G$ be a finite simplicial directed graph. The
admissible boundary $\partial _{a}$ of $G$ is defined by a boundary with the
following condition;

\strut

\begin{center}
$\{v_{1},v_{2}\}\subseteq \partial _{a}\Longleftrightarrow $ there is no
finite path connecting $v_{1}$ and $v_{2}$ on $G.$
\end{center}

\strut

The admissible boundary quotient $\partial _{a}$ is in general not an
invariant on finite simplicial directed graphs. However, we show that the
admissible boundary quotient is an invariant on the finite simplicial 
\textbf{connected} directed graph \textbf{with mixed maximal type}.

\strut

In Chapter 2, we introduce the subgraph boundaries which is different from
the boundaries defined in Chapter 1. Let $G$ be a finite directed graph (not
necesarily simplicial) and let $H$ be a full subgraph of $G.$ Then the
boundary $\partial _{H}$ $=$ $V(H)$ $\cup $ $E(H)$ is called the subgraph
boundary of $G$ with respect to $H.$ Similar to Chapter 1, by identifying
all elements in $\partial _{H}$ to the base point $v_{\partial _{H}},$ we
can construct the subgraph boundary quotient graph $G$ $/$ $\partial _{H}.$
Notice that,

\strut

\begin{center}
$\left| V(G\text{ }/\text{ }\partial _{H})\right| =\left| V(G)\right|
-\left| V(H)\right| +1,$

$\left| E(G\text{ }/\text{ }\partial _{H})\right| =\left| E(G)\right|
-\left| E(H)\right| $
\end{center}

and hence

\begin{center}
$\left| V(G)\cup E(G)\right| =\left| V(G)\cup E(G)\right| -\left| V(H)\cup
E(H)\right| +1.$
\end{center}

\strut

By using the above quantity, we will define the subgraph boundary index $%
Ind_{G}(H)$ of a full subgraph $H$ of $G,$ by

\strut

\begin{center}
$Ind_{G}(H)=\exp \left( \left| G\text{ }/\text{ }\partial _{H}\right|
-1\right) .$
\end{center}

\strut

If we denote $1$ as the trivial graph (the one-vertex-no-edge graph), then
we can get that

\strut

\begin{center}
$Ind_{G}(H)=\frac{Ind_{G}(1)}{Ind_{H}(1)}.$
\end{center}

\strut

We will consider this quantity in detail in Chapter 2.

\strut \strut

\strut \strut

\strut \strut

\section{Boundary Quotient Graphs}

\strut

\strut

In this chapter, we will define a boundary quotient graph $G$ $/$ $\partial $
and observe some properties of it. Throughout this chapter, let $G$ be a
finite simplicial directed graph with its vertex set $V(G)$ and its edge set 
$E(G).$ Since the graph $G$ is simplicial, this graph contains neither
multiple edges between two vertices nor loop-edges. For example, all growing
trees and all one-flow circulant graphs are simplicial graphs. The boundary $%
\partial $ of $G$ is defined by a subset of $V(G).$ The boundary quotient
graph $G$ $/$ $\partial $ is defined by a graph $Q$ under the boundary
quotient, also denoted by $\partial .$

\strut

\strut

\strut

\subsection{Boundary Quotient Graphs}

\strut

\strut \strut

Let $G$ be the given finite simplicial directed graph.

\strut

\begin{definition}
A boundary $\partial $ of $G$ is a subset of $V(G).$ If $\partial $ $=$ $%
V(G),$ we will say that $\partial $ is total in $G$. Otherwise, we say that $%
\partial $ is proper in $G.$ The boundary quotient of $\partial $ is defined
by the following relations;

\strut

(R1) All vertices in $\partial $ are identified to the point $v_{\partial }.$
This new point $v_{\partial }$ is called the base point of $\partial .$
i.e., if $v_{1}$ $\neq $ $v_{2}$ in $\partial ,$ then identify $v_{1}$ and $%
v_{2}$ with the base point $v_{\partial }.$

\strut

(R2) If $v_{1}$ $\neq $ $v_{2}$ in $\partial $ and if there exists a direct
edge $e$ connecting $v_{1}$ and $v_{2},$ then identify $e$ to the loop-edge $%
l_{e}$ at the base point $v_{\partial }.$

\strut

The boundary quotient with (R1) and (R2) is also denoted by $\partial ,$ if
there is no confusion.
\end{definition}

\strut

We say that the edges of $G,$ which are not affected by the boundary
quotient $\partial ,$ are non-$\partial $ edges. Also, the loop-edges
constructed by (R2) of $\partial ,$ are said to be $\partial $\textbf{-}%
loop-edges at the base point\textbf{\ }$v_{\partial }.$ With respect to the
fixed boundary quotient $\partial $ on the given graph $G,$ we will define
the corresponding boundary quotient graph $G$ $/$ $\partial .$

\strut

\begin{definition}
The boundary quotient graph $G$ $/$ $\partial $ is a finite directed
(non-simplicial) graph $Q$ with its vertex set

\strut

\ $\ \ \ \ \ \ \ \ \ \ \ \ \ \ \ V(Q)=\{v_{\partial }\}\cup \{v\in V(G):v\in
V(G)$ $\setminus $ $\partial \}$

\strut

and its edge set

\strut

$\ \ \ \ \ \ \ E(Q)=\{e\in E(G):e$ is non-$\partial \}\cup \{l:l$ is a $%
\partial $-loop-edge$\}.$
\end{definition}

\strut

From now, for convenience, we will denote $e$ $=$ $v_{1}ev_{2},$ where $e$
is an edge connecting $v_{1}$ and $v_{2}$ with the direction from $v_{1}$ to 
$v_{2}.$ (i.e., $e$ is an edge with its initial vertex $v_{1}$ and its
terminal vertex $v_{2}.$) Also, if $w$ $=$ $e_{1}$ ... $e_{k}$ is a finite
path with its length $k$, with the admissible edges $e_{1},$ ..., $e_{k}$
and if $e_{1}$ $=$ $v_{1}e_{1}v^{\prime }$ and $e_{k}$ $=$ $v^{\prime \prime
}e_{k}v_{2},$ then we write $w$ $=$ $v_{1}wv_{2}$ to emphasize the initial
and terminal vertices of $w.$ We give some fundamental examples of boundary
quotient graphs.

\strut \strut

\begin{example}
Let $G$ be a graph with $V(G)$ $=$ $\{v_{1},$ $v_{2}\}$ and $E(G)$ $=$ $\{e$ 
$=$ $v_{1}ev_{2}\}.$ Suppose we have the total boundary $\partial $ $=$ $%
V(G).$ Then, by the boundary quotient $\partial ,$ we can get the boundary
quotient graph $G$ $/$ $\partial ,$ as a graph $Q$ with its vertex set $V(Q)$
$=$ $\{v_{\partial }\}$ and $E(Q)$ $=$ $\{l_{e}$ $=$ $v_{\partial }$ $l_{e}$ 
$v_{\partial }\}.$ i.e., the graph $Q$ is the one-vertex-one-loop-edge
graph. The $\partial $-loop edge $l_{e}$ is constructed from the edge $e.$
\end{example}

\strut

\begin{example}
Let $G$ be a graph with $V(G)$ $=$ $\{v_{1},$ $v_{2},$ $v_{3}\}$ and $E(G)$ $%
=$ $\{e_{1}$ $=$ $v_{1}$ $e_{1}$ $v_{2},$ $e_{2}$ $=$ $v_{1}$ $e_{2}$ $%
v_{3}\}.$ Suppose we have a $\partial _{12}$ $=$ $\{v_{1},$ $v_{2}\}.$ Then
the boundary quotient graph $G$ $/$ $\partial _{12}$ is a graph $Q_{12}$
with $V(Q_{12})$ $=$ $\{v_{\partial _{12}},$ $v_{3}\}$ and $E(Q_{12})$ $=$ $%
\{l_{e_{1}}$ $=$ $v_{\partial _{12}}$ $l_{e_{1}}$ $v_{\partial _{12}},$ $%
e_{2}\}.$ Similarly, if we have the boundary $\partial _{13}$ $=$ $\{v_{1},$ 
$v_{3}\},$ then we have the boundary quotient graph $G$ $/$ $\partial _{13}$%
, as a graph $Q_{13}$ with $V(Q_{13})$ $=$ $\{v_{\partial _{13}},$ $v_{2}\}$
and $E(Q_{13})$ $=$ $\{l_{e_{2}}$ $=$ $v_{\partial _{13}}$ $l_{e_{2}}$ $%
v_{\partial _{13}},$ $e_{1}\}.$ Now, the boundary $\partial $ is total in $%
G. $ i.e., $\partial $ $=$ $V(G).$ Then the boundary quotient graph $G$ $/$ $%
\partial $ is the graph $Q$ with $V(Q)$ $=$ $\{v_{\partial }\}$ and $E(Q)$ $%
= $ $\{l_{e_{1}}$ $=$ $v_{\partial }$ $l_{e_{1}}$ $v_{\partial },$ $\
l_{e_{2}} $ $=$ $v_{\partial }$ $l_{e_{2}}$ $v_{\partial }\}.$
\end{example}

\strut

Define the one-flow circulant graph $C_{N}$, as a graph $K$ with its vertex
set $V(K)$ $=$ $\{v_{1},$ ..., $v_{N}\}$ and its edge set

\strut

\begin{center}
$E(K)=\{e_{j}=v_{j}e_{j}v_{j+1}:j=1,...,N,$ and $v_{N+1}\overset{def}{=}%
v_{1}\}.$
\end{center}

\strut

\begin{example}
Let $G$ be a one-flow circulant graph $C_{3}.$ Suppose we have the total
boundary $\partial $ $=$ $V(G).$ Then the corresponding boundary quotient
graph $G$ $/$ $\partial $ is the graph $Q$ with $V(G)$ $=$ $\{v_{\partial
}\} $ and $E(G)$ $=$ $\{l_{e_{j}}$ $=$ $v_{\partial }$ $l_{e_{j}}$ $%
v_{\partial } $ $:$ $j$ $=$ $1,$ $2,$ $3\}.$ So, the graph $G$ $/$ $\partial 
$ is the one-vertex-three-loop-edge graph.
\end{example}

\strut \strut

Let $G_{1}$ and $G_{2}$ be finite simplicial directed graphs and let $v_{1}$
and $v_{2}$ be arbitrary fixed vertices of $G_{1}$ and $G_{2},$
respectively. Define the \textbf{vertex-fixed} \textbf{glued graph} $G_{1}$ $%
\#_{v_{1}\#v_{2}}$ $G_{2}$ \textbf{at the glued vertex} $v_{1}$ $\#$ $v_{2}$
by the graph $G$ with the following conditions;

\strut

(C1) \ Identify $v_{1}$ and $v_{2}.$ This identified vertex $v_{1}\#v_{2}$
in $G$ is called the glued

\ \ \ \ \ \ \ \ vertex.

(C2) $\ V(G)=\{v_{1}\#v_{2}\}\cup \left( V(G_{1})\text{ }\setminus \text{ }%
\{v_{1}\}\right) \cup \left( V(G_{2})\text{ }\setminus \text{ }%
\{v_{2}\}\right) .$

(C3) $\ E(G)=E(G_{1})\cup E(G_{2}).$\strut

\strut

Inductively, we can have the vertex-fixed glued graph $G_{1}$ $%
\#_{v_{1}\#...\#v_{m}}$ ... $\#_{v_{1}\#...\#v_{m}}$ $G_{m},$ for $m$ $\in $ 
$\Bbb{N},$ at its glued vertex $v_{1}$ $\#$ ... $\#$ $v_{m}.$ The
vertex-fixed glued graphs are depending on the choice of their glued
vertices. i.e., the vertex-fixed glued graph $G_{1}$ $\#_{v}$ ... $\#_{v}$ $%
G_{m}$ at the glued vertex $v$ $=$ $v_{1}$ $\#$ ... $\#$ $v_{m}$ and the
vertex-fixed glued graph $G_{1}$ $\#_{v^{\prime }}$ ... $\#_{v^{\prime }}$ $%
G_{m}$ at the glued vertex $v^{\prime }$ $=$ $v_{1}^{\prime }$ $\#$ ... $\#$ 
$v_{m}^{\prime }$ are not graph-isomorphic, in general, if $v_{j}$ $\neq $ $%
v_{j}^{\prime }$, for some $j$ in $\{1,$ ..., $m\}.$

\strut

\begin{proposition}
Let $G_{j}$ be one-flow circulant graph $C_{n_{j}},$ for all $j$ $=$ $1,$
..., $m.$ Then the vertex-fixed glued graph $G_{1}$ $\#$ ... $\#$ $G_{m}$ is
independent of the choice of glued vertex, up to graph-isomorphisms. $%
\square $
\end{proposition}

\strut

The above proposition is easily proved, by the definition of the circulant
graphs and graph-isomorphisms.

\strut

\begin{example}
Let $G$ be a one-flow circulant graph $C_{4}$ and let the boundary $\partial 
$ is $\{v_{1},$ $v_{2}\}.$ Then the boundary quotient graph $G$ $/$ $%
\partial $ is the graph $Q$ with $V(Q)$ $=$ $\{v_{\partial },$ $v_{3},$ $%
v_{4}\}$ and $E(G)$ $=$ $\{l_{e_{1}}$ $=$ $v_{\partial }$ $l_{e_{1}}$ $%
v_{\partial },$ $e_{2},$ $e_{3},$ $e_{4}\}.$ Notice that this graph $Q$ is
the vertex-fixed glued graph $G_{1}$ $\#$ $G_{2},$ where $G_{1}$ is the
graph with $V(G)$ $=$ $\{v_{\partial }\}$ and $E(G)$ $=$ $\{l_{e_{1}}\}$ and 
$G_{2}$ is the graph with $V(G_{2})$ $=$ $\{v_{\partial },$ $v_{2},$ $%
v_{3}\} $ and $E(G)$ $=$ $\{e_{2},$ $e_{3},$ $e_{4}\}.$ Moreover, the graph $%
G_{1}$ is graph-isomorphic to the one-vertex-one-edge graph and the graph $%
G_{2}$ is graph-isomorphic to the one-flow circulant graph $C_{3}.$
\end{example}

\strut \strut

We will define the following special types on the finite directed simplicial
graphs.

\strut

\begin{definition}
(1) Let $G$ be a graph which is graph-isomorphic to $C_{n}$, the one-flow
circulant graphs with $n$-vertices. Then the graph $G$ is said to be of type 
$C_{n}.$

\strut

(2) We say that a graph $G$ is of type $L_{m},$ if the graph $G$ is
graph-isomorphic to the one-vertex-$m$-loop-edge graph.

\strut

(3) A graph $G$ is of type $T$ if $G$ is isomorphic to a directed tree $T.$

\strut

(4) A graph $G$ is of mixed type if there exists full subgraphs $G_{1},$
..., $G_{n}$ of $G$ such that (i) $\{G_{1},$ ..., $G_{n}\}$ is the minimal
covering of $G$ and (ii) $G_{j}$ is of type $C_{N}$ or of type $L_{m}$ or of
type $T$, for $j$ $=$ $1,$ ..., $n.$ (Since a graph $G$ is finite, we can
always choose such finite covering consisted of full subgraphs.)
\end{definition}

\strut \strut

\begin{proposition}
Let $G$ be a one-flow circulant graph $C_{N}$ and let $\partial $ be the
total boundary of $G.$ Then the boundary quotient graph $G$ $/$ $\partial $
is the graph of type $L_{N}.$
\end{proposition}

\strut

\begin{proof}
Let $G$ be a one-flow circulant graph and let $\partial $ $=$ $V(G)$ be the
total boundary of $G.$ Then, by the boundary quotient, the graph $G$ $/$ $%
\partial $ is a graph $Q$ with $V(Q)$ $=$ $\{v_{\partial }\}$ and $E(Q)$ $=$ 
$\{l_{e_{j}}$ $=$ $v_{\partial }$ $l_{e_{j}}$ $v_{\partial }$ $:$ $j$ $=$ $%
1, $ ..., $N\},$ where $e_{j}$'s are edges $e_{j}$ $=$ $v_{j}$ $e_{j}$ $%
v_{j+1}$ in $E(G),$ for all $j$ $=$ $1,$ ..., $N,$ with $v_{N+1}$ $\overset{%
let}{=}$ $v_{1}.$ Therefore, the graph $G$ $/$ $\partial $ is of type $L_{N}.
$
\end{proof}

\strut

More generally, we can get that;

\strut

\begin{proposition}
(1) Let $G$ be a one-flow circulant graph $C_{N}$ and let $\partial _{k}$ $=$
$\{v_{j},$ $v_{j+1},$ ..., $v_{j+k-1}\}$ be a boundary consisting of $k$%
-vertices, where $1$ $\leq $ $j$ $<$ $N$ and $2$ $\leq $ $j$ $+$ $k$ $<$ $N$ 
$+$ $1.$ Then the boundary quotient graph $G$ $/$ $\partial $ is of mixed
type $(L_{k}$ , \ $C_{N-k}).$

\strut \strut \strut

(2) Let $G$ be a one-flow circulant graph $C_{N}$ and let $\partial _{k}$ $=$
$\{v_{j_{1}},$ $v_{j_{2}}$..., $v_{j_{n}}\}$, where $(j_{1},$ ..., $j_{n})$
is a sequence in $\{1,$ ..., $N\},$ for $n$ $<$ $N$ satisfying that $j_{k+1}$
$=$ $j_{k}$ $+$ $t_{k},$ for $t_{k}$ $>$ $1$ and for $k$ $=$ $1,$ ..., $n.$
Then the boundary quotient graph $G$ $/$ $\partial _{k}$ is of type $%
(C_{j_{1}-1},$ $C_{t_{2}-1},$ $C_{t_{3}-1},$ ..., $C_{N-j_{n}-1},$ $%
C_{N-j_{n}}).$
\end{proposition}

\strut

\begin{proof}
(1) Suppose $1$ $<$ $j$ $<$ $N$ and $1$ $<$ $j+k$ $-$ $1$ $<$ $N.$ Consider
the full subgraph $L$ of $G$ with its vertex set $V(L)$ $=$ $\partial _{k}.$
Then, by the previous proposition, the boundary quotient graph $L$ $/$ $%
\partial $ is of type $L_{k}.$ With respect to the base point $v_{\partial
}, $ there exists subgraph $K$ of type $(N$ $-k)$ such that $G$ $/$ $%
\partial $ $=$ $L$ $\#_{v_{\partial }}$ $K.$ i.e., the boundary quotient
graph $G$ $/$ $\partial $ is the glued graph $L$ $\#_{v_{\partial }}$ $K$ of
maximal type, with the base point $v_{\partial },$ where $G_{1}$ is of type $%
L_{k}$ and $G_{2}$ is of type $C_{N-k}.$

\strut

(2) The boundary quotient graph $G$ $/$ $\partial _{k}$ satisfies that $G$ $%
/ $ $\partial _{k}$ is graph-isomorphic to the glued graph $K_{1}$ $%
\#_{v_{\partial }}$ $K_{2}$ $\#_{v_{\partial }}$ ... $\#_{v_{\partial }}$ $%
K_{n},$ with its glued vertex (which is the base point of $G$ $/$ $\partial
_{k}$), where $K_{1}$ is of type $C_{j_{1}-1}$ and $K_{i}$ is of type $%
C_{t_{i}-1},$ for each $i$ $=$ $1,$ ..., $n$ $-$ $1,$ and $K_{n}$ is of type 
$C_{N-j_{n}}.$
\end{proof}

\strut

Let $T_{N}$ be a finite directed tree with $N$-vertices. The given tree $%
T_{N}$ is said to be a growing tree if there exists the root vertex $v_{0},$
and there always exists a unique finite path $w_{v}$ $=$ $v_{0}$ $w$ $v$ on $%
T_{N},$ for all other vertex $v$ $\neq $ $v_{0}$ in $V(G).$ Thus the growing
tree $T_{N}$ has only one-flow from the root vertex $v_{0}$.

\strut \strut \strut \strut

\begin{proposition}
Let $T_{N}$ be a tree with $N$-vertices and let $\partial $ be the total
boundary of $T_{N}.$ Then the boundary quotient graph $T_{N}$ $/$ $\partial $
is of type $L_{N-1}.$ $\square $
\end{proposition}

\strut

The above proposition is proved by induction. The above proposition says
that if $T$ is a tree and if $\partial $ $=$ $V(T)$ is the total boundary,
then $T$ is of type $L_{\left| E(T)\right| }.$ Also, it says that the total
boundary quotient is not an invariants on finite simplicial graphs. Since
the total boundary quotient of $C_{N-1}$ is also type $L_{N-1}.$ Trivially
the graphs $C_{N-1}$ and $T_{N}$ are not graph-isomorphic. But their total
boundary quotient graphs are of type $L_{N-1}.$ Thus the total boundary
quotient is not an invariants on graphs.

\strut \strut \strut \strut

Notice that every finite connected simplicial directed graph $G$ is
graph-isomorphic the \textbf{iterated glued graph}

\strut

\begin{center}
$G_{1}\#_{v_{1}}\left( G_{2}\#_{v_{2}}\left(
G_{3}\#_{v_{3}}(G_{4}\#_{v_{4}}...(G_{n}\#_{v_{n}}G_{n+1})\right) \right) $
\end{center}

\strut

of trees, circulant graphs and trivial graphs. Different from the
vertex-fixed glued graph case, it is not necessary that $v_{1}$ $=$ $v_{2}$ $%
...$ $=$ $v_{n}$. Remark that the iterated glued graphs of $G$ are not
uniquely determined. But we can choose the maximal one-flow circulant full
subgraphs in $G$ and the maximal sub-trees in $G,$ in their neighborhoods.
Also, if we fix the subset $V$ of the vertex set $V(G),$ then we can choose
the maximal one-flow circulant full subgraphs and the maximal full sub-trees
in $G,$ with respect to the fixed vertices in $V$. For example, let $G$ be a
graph with

\strut

\begin{center}
$V(G)$ $=$ $\{v_{1},$ ..., $v_{8}\}$
\end{center}

and

\begin{center}
$E(G)=\left\{ 
\begin{array}{c}
e_{1}=v_{1}e_{1}v_{2},\text{ \ }e_{2}=v_{2}e_{2}v_{3}, \\ 
e_{3}=v_{3}e_{3}v_{5},\text{ \ }e_{4}=v_{5}e_{4}v_{4}, \\ 
e_{5}=v_{4}e_{5}v_{1},\text{ \ }e_{6}=v_{5}e_{6}v_{6}, \\ 
e_{7}=v_{6}e_{7}v_{7},\text{ \ }e_{8}=v_{6}e_{8}v_{8}
\end{array}
\right\} .$
\end{center}

\strut

Fix a boundary $\partial _{1}$ $=$ $\{v_{2},$ $v_{5},$ $v_{8}\}.$ Then, for
the fixed vertices $v_{2},$ $v_{5}$ and $v_{8},$ we can choose the full
subgraphs $G_{1},$ $G_{2}$ and $G_{3},$ as graphs with

\strut

$\ \ \ \ \ \ \ \ \ V(G_{1})=\{v_{2},v_{3},v_{5}\}$ \ and \ $%
E(G_{1})=\{e_{2},e_{3}\},$

$\ \ \ \ \ \ \ \ \ V(G_{2})=\{v_{1},v_{2},v_{4},v_{5}\}$ \ and \ $%
E(G_{2})=\{e_{1},e_{4},e_{5}\},$

and

$\ \ \ \ \ \ \ \ \ V(G_{3})=\{v_{5},v_{6},v_{7},v_{8}\}$ \ and \ $%
E(G_{3})=\{e_{6},e_{7},e_{8}\}.$

\strut

All those three full subgraphs are trees. Then the boundary quotient graph $%
G $ $/$ $\partial _{1}$ is the glued graph

\strut

\begin{center}
$G_{1}^{\prime }\ \#_{v_{\partial _{1}}}\ \left( G_{2}^{\prime }\
\#_{v_{\partial _{1}}}\ \left( G_{3}^{\prime }\#_{v_{6}}G_{4}\right) \right)
,$
\end{center}

\strut

where $G_{1}^{\prime }$ is of type $C_{2},$ $G_{2}^{\prime }$ is of type $%
C_{3}$ and $G_{3}^{\prime }$ is of type $C_{2}$ and where $G_{4}$ is a tree
with $V(G_{4})$ $=$ $\{v_{6},$ $v_{7}\}$ and $E(G_{4})$ $=$ $\{e_{7}\}.$
Also, we have the following iterated glued graph of the above given graph $G$%
;

\strut

\begin{center}
$K_{1}\#_{v_{2}}\left( K_{2}\#_{v_{3}}\left( \left( \left(
K_{3}\#_{v_{4}}\left( K_{4}\#_{v_{5}}K_{5}\#_{v_{1}}K_{6}\right) \right)
\#_{v_{4}}K_{7}\right) \#_{v_{6}}K_{8}\right) \right) ,$
\end{center}

\strut

where

\strut

\begin{center}
$V(K_{1})=\{v_{1},v_{2}\}$ \ and \ $E(K_{1})=\{e_{1}\},$

$V(K_{2})=\{v_{2},v_{3}\}$ \ and \ $E(K_{2})=\{e_{2}\},$

$V(K_{3})=\{v_{3},v_{4}\}$ \ and \ $E(K_{3})=\{e_{3}\},$

$V(K_{4})=\{v_{4},v_{5},v_{6}\}$ and $E(K_{4})=\{e_{4},e_{6}\},$

$V(K_{5})=\{v_{1},v_{5}\}$ \ and \ $E(K_{5})=\{e_{5}\},$

$V(K_{6})=\{v_{1}\}$ \ \ \ \ \ \ and \ \ \ \ $\ E(K_{6})=\emptyset ,$

$V(K_{7})=\{v_{4},v_{6}\}$ \ and \ $E(K_{7})=\{e^{6}\}$
\end{center}

and

\begin{center}
$V(K_{8})=\{v_{6},v_{7},v_{8}\}$ and $E(K_{8})=\{e_{7},e_{8}\}.$
\end{center}

\strut

The above example shows us that the iterated glued graph of the
graph-isomorphic graph is not uniquely determined.. However, we can choose
the iterated glued graph of a given graph by the maximal $CT$-iterated glued
graph.

\strut

\begin{definition}
Let $G$ be a finite connected simplicial directed graph. A $CT$-iterated
glued graph of $G$ is the graph-isomorphic graph $K$ $=$ $K_{1}$ $\#_{v_{1}}$
$K_{2}$ ... $\#_{v_{n}}K_{n+1},$ where each $K_{j}$ is a full subgraph of $G$
which is of type $C_{N}$ or of type $T$ or of type $1$, for $j$ $=$ $1,$
..., $n$ $+$ $1.$ A graph is of type $1,$ if it is the trivial graph which
is the one-vertex-no-edge graph. We say that the $CT$-iterated glued graph $%
K $ is maximal if each gluing components $K_{j}$'s are the maximal full
subgraphs satisfying the type, in its neighborhood of $G.$ The maximal type $%
(t_{1},$ ..., $t_{n+1})$ of $G$ is defined by

\strut

$\ \ \ \ \ \ \ \ \ \ \ \ \ \ \ t_{j}=\left\{ 
\begin{array}{lll}
C_{N_{j}} &  & \text{if }K_{j}\text{ is of type }C_{N_{j}} \\ 
T &  & \text{if }K_{j}\text{ is of type }T \\ 
1 &  & \text{if }K_{j}\text{ is the trivial graph,}
\end{array}
\right. $

\strut

for $j$ $=$ $1,$ ..., $n+1.$ (Recall that, by definition, $C_{N}$ is the
one-flow circulant graph and $T$ is a tree.) The type of the maximal $CT$%
-iterated glued graph is called the maximal type of $G.$
\end{definition}

\strut \strut

In the above example, the first iterated glued graph of $G$ is a $CT$%
-iterated glued graph which is not maximal. The second iterated graph of $G$
is an iterated glued graph of $G$ which is not a $CT$-iterated glued graph.
The following glued graph of $G$ is the maximal $CT$-iterated glued graph of 
$G$;

\strut

\begin{center}
$G_{1}$ $\ \#_{v_{5}}$ $\ G_{2}$
\end{center}

with

\begin{center}
$V(G_{1})=\{v_{1},v_{2},v_{3},v_{4},v_{5}\}$ and $E(G_{1})=%
\{e_{1},e_{2},e_{3},e_{4},e_{5}\}$
\end{center}

and

\begin{center}
$V(G_{2})=\{v_{5},v_{6},v_{7},v_{8}\}$ \ \ \ \ \ \ and \ \ \ \ \ $%
E(G_{2})=\{e_{6},e_{7},e_{8}\},$
\end{center}

\strut

where $G_{1}$ is of type $C_{5}$ and $G_{2}$ is of type $T_{4},$ where $%
T_{4} $ is the graph with $V(T_{4})$ $=$ $\{b_{1},$ $b_{2},$ $b_{3},$ $%
b_{4}\}$ and $E(T_{4})$ $=$ $\{f_{1},$ $f_{2},$ $f_{3}\},$ where $f_{1}$ $=$ 
$b_{1}$ $f_{1}$ $b_{2},$ $f_{2}$ $=$ $b_{2}$ $f_{2}$ $b_{3}$ and $f_{3}$ $=$ 
$b_{2}$ $f_{3}$ $b_{4}.$ i.e., the maximal $CT$-iterated glued graph $G_{1}$ 
$\#_{v_{5}}$ $G_{2}$ of $G$ is of its maximal type $(C_{5},$ $T_{4}).$

\strut

In the above definition, we introduced the trivial graph. The trivial graphs
are needed for finding the maximal $CT$-iterated glued graph. For example,
let $G$ be a graph with $V(G)$ $=$ $\{v_{1},$ $v_{2}\}$ and $E(G)$ $=$ $%
\{e_{1}$ $=$ $v_{1}$ $e_{1}$ $v_{2},$ $e_{2}$ $=$ $v_{2}$ $e_{2}$ $v_{1}\}.$
This graph has the following maximal $CT$-iterated glued graph $K,$

\strut

\begin{center}
$K$ $=$ $K_{1}$ $\#_{v_{2}}$ $K_{2}$ $\#_{v_{1}}$ $K_{3},$
\end{center}

with

\begin{center}
$V(K_{1})$ $=$ $\{v_{1},$ $v_{2}\}$ \ and \ $E(K_{1})=\{e_{1}\},$

$V(K_{2})=\{v_{1},v_{2}\}$ \ and \ $E(K_{2})=\{e_{2}\}$
\end{center}

and

\begin{center}
$V(K_{3})=\{v_{1}\}$ \ \ and \ \ $E(K_{3})=\varnothing .$
\end{center}

\strut

Remark that $K_{3}$ is the trivial graph of its type $1.$ \strut So, this
graph $G$ is of mixed maximal type $(T,$ $T^{\prime },$ $1),$ where $T$ $=$ $%
K_{1}$ and $T^{\prime }$ $=$ $K_{2}.$ We can have the following lemma.

\strut

\begin{lemma}
Let $G$ be a finite simplicial directed graph. Then there always exists the
maximal $CT$-iterated glued graph. Futhermore, the maximal $CT$-iterated
glued graph is uniquely determined, up to graph-isomorphisms.
\end{lemma}

\strut

\begin{proof}
Let $G$ be the given graph. Then this graph is the disjoint union of $G_{1},$
..., $G_{n},$ where $G_{j}$ is the connected component of $G,$ for $j$ $=$ $%
1,$ ..., $n.$ Since $G$ is finite simplicial, all $G_{j}$'s are also finite
simplicial. By the simplicity of $G_{j},$ this graph $G_{j}$ does not
contain the loop-edges and multiple edges. So, graphically, we can choose
the full subgraphs $G_{j1},$ ..., $G_{jn_{j}},$ where $G_{ji}$ is of type $%
C_{N_{i}}$ or of type $T_{k_{i}}$ or og type $1,$ for $i$ $=$ $1,$ ..., $%
n_{j}$ and for all $j$ $=$ $1,$ ..., $n.$ By the Axiom of Choice, we can
choose such maximal family.

\strut

Now, assume that $K_{1}$ and $K_{2}$ are maximal $CT$-iterated glued graphs
of $G.$ For convenience, suppose that the given graph $G$ is connected.
Let's assume that $K_{1}$ and $K_{2}$ are not graph-isomorphic. Without loss
of generality, we may suppose that $K_{1}$ (resp. $K_{2}$) has $K_{i_{1}},$
..., $K_{i_{s}}$ (resp. $K_{j_{1}},$ ..., $K_{j_{t}}$) gluing components
which are of type $C_{N_{i_{1}}},$ ..., $C_{N_{i_{s}}}$ (resp. $%
C_{N_{j_{1}}},$ ..., $C_{N_{j_{t}}}$), respectively, and $s$ $\neq $ $t$ in $%
\Bbb{N}.$ Also, assume that $s$ $<$ $t.$ By the maximality, $K_{2}$ cannot
be the maximal $CT$-iterated glued graph of $G.$
\end{proof}

\strut \strut

The above lemma says that every finite simplicial directed graph is
graph-isomorphic to the maximal $CT$-iterated glued graph of maximal type $%
((C_{1},$ ..., $C_{n}),$ $(T_{1},$ ..., $T_{k}),$ $(1,$ ..., $1)).$
Therefore, we can conclude that the boundary quotient graph $G$ $/$ $%
\partial $ is, in general, of maximal type $((L_{1},$ ..., $L_{m}),$ $%
(C_{1}, $ ..., $C_{n}),$ $(T_{1},$ ..., $T_{k}),$ $(1,$ ..., $1)).$

\strut

\begin{theorem}
Let $G$ be an arbitrary finite simplicial graph and let $\partial $ be a
boundary of $G.$ Then, in general, the quotient graph $G$ $/$ $\partial $ is
of maximal type $((L_{m_{1}},$ ... $L_{m_{s}}),$ $(C_{N_{1}},$ ..., $%
C_{N_{n}}),$ $(T_{k_{1}},$ ..., $T_{k_{m}}),$ $(1,$ ..., $1)),$ where $%
\left| E(G)\right| $ $=$ $\sum_{p=1}^{s}$ $m_{p}$ $+$ $\sum_{i=1}^{n}$ $%
N_{i} $ $+$ $\sum_{r=1}^{m}$ $k_{r}$.
\end{theorem}

\strut

\begin{proof}
Let $G$ and $\partial $ be given and assume that $G$ is connected. By the
previous lemma, the graph $G$ is graph-isomorphic to the maximal $CT$%
-iterated glued graph of $G_{1},$ ..., $G_{n},$ where $G_{j}$ is either of
type $C_{n_{j}}$ or of type $T_{k_{j}},$ for $j$ $=$ $1,$ ..., $n.$ By
regarding each gluing component $G_{j}$ as the full subgraph $K_{j}$ of $G,$
we can get the subboundary $\partial _{j},$ for each $j.$ i.e., $\partial
_{j}$ $=$ $V(G_{j})$ $\cap $ $\partial .$ Notice that we have the boundary
quotient graph $G$ $/$ $\partial $ is the vertex-fixed glued graph of $G$ $/$
$\partial _{j},$ $j$ $=$ $1,$ ..., $n,$ with its glued vertex $v_{\partial }$
$=$ $v_{\partial _{1}}$ $\#$ ... $\#$ $v_{\partial _{n}}.$ By the previous
propositions, if $G_{j}$ is of type $C_{N_{j}}$ (or of type $T_{k_{j}}$),
then $G_{j}$ $/$ $\partial _{j}$ is of maximal type of $(L_{m_{j}},$ $%
(C_{N_{j},1},$ ..., $C_{N_{j},n_{j}})),$ (resp. of maximal type $(L_{m_{j}},$
$(C_{N_{j},1},$ ..., $C_{N_{j},m_{j}}),$ $(T_{j,k_{1}},$ ..., $T_{j,k_{j}}))$%
) for $j$ $=$ $1,$ ..., $n.$ Since $G$ $/$ $\partial $ $=$ $(G_{1}$ $/$ $%
\partial _{1})$ $\#_{v_{\partial }}$ ... $\#_{v_{\partial }}$ $(G_{n}$ $/$ $%
\partial _{n}),$ it is of maximal type with $L_{m}$'s, $C_{N}$'s and $T$'s.
Remark that $\left| E(G\text{ }/\text{ }\partial )\right| $ $=$ $\left|
E(G)\right| $. Therefore,

\strut

$\ \ \ \ \ \ \ \ \ \ \ \ \ \ \ \ \ \ \ \left| E(G)\right| $ $=$ $m$ $+$ $%
\sum_{i=1}^{n}$ $N_{i}$ $+$ $\sum_{r=1}^{m}$ $k_{r}.$

\strut

Now, suppose that the graph $G$ is the disjoint union of finite simplicial
connected directed graphs $G_{1},$ ..., $G_{k}.$ Then, for each $G_{k},$ we
have the above result. And hence we can get the desired result.
\end{proof}

\strut \strut

If we consider the total boundary quotient graph, we have the following
simple results;

\strut

\begin{theorem}
Let $G$ be a finite connected simplicial directed graph and let $\partial $
be the total boundary of $G.$ Then the boundary quotient graph $G$ $/$ $%
\partial $ is of type $L_{\left| E(G)\right| }.$
\end{theorem}

\strut

\begin{proof}
By the previous lemma, $G$ is graph-isomorphic to the maximal $CT$-iterated
glued graph $K$ with its gluing components $K_{1},$ ..., $K_{n}.$ Then each $%
K_{j}$ is of type $C_{N}$ or type $T$ or of type $1.$ Notice that if $%
\partial $ is total in $G,$ then

\strut

$\ \ \ \ \ \ \ \ \ \ \ \ \ \ \ \ \ G$ $/$ $\partial =\left( K_{1}\text{ }/%
\text{ }\partial _{1}\right) \#$ $...$ $\#$ $\left( K_{n}\text{ }/\text{ }%
\partial _{n}\right) ,$

\strut

where the right-hand side is the vertex-fixed glued graph and where $%
\partial _{j}$ $=$ $\partial \cap V(K_{j})$ $=$ $V(K_{j}),$ for $j$ $=$ $1,$
..., $n.$ If $K_{j}$ is trivial, then $K_{j}$ $/$ $\partial _{j}$ is
trivial. If $K_{j}$ is either of type $C_{N}$ or of type $T,$ then $K_{j}$ $%
/ $ $\partial _{j}$ is of type $L_{\left| E(K_{j})\right| }.$ Therefore, we
can conclude the result.
\end{proof}

\strut \strut \strut \strut

\strut \strut \strut

\strut

\subsection{General Total Boundary Quotient Graphs}

\strut

\strut

In the previous section, we only considered the finite simplicial directed
graphs. We say that the directed graph $G$ is finite if $\left| V(G)\right| $
$<$ $\infty $ and $\left| E(G)\right| $ $<$ $\infty .$ However, they may
have the loop-edges and multiple edges between two vertices. However, we can
extend the above results in Section 1.1 to the general finite graph cases.
Moreover, we have that

\strut \strut \strut

\begin{theorem}
Let $G$ be a finite directed graph and $\partial $, the total boundary of $G.
$ Then the total boundary quotient graph $G$ $/$ $\partial $ of $G$ is of
type $L_{\left| E(G)\right| }.$ $\square $
\end{theorem}

\strut \strut

We cannot use the maximal $CT$-iterated glued graph technique to prove the
above geneal case. But we can use the edge-iterated glued graph of $G.$ If
there are multiple edges $e_{1},$ ..., $e_{k}$ connecting $v_{1}$ and $v_{2},
$ with same direction. Take the full subgraph $K_{j}$ with $V(K_{j})$ $=$ $%
\{v_{1},$ $v_{2}\}$ and $E(K_{j})$ $=$ $\{e_{j}\},$ for $j$ $=$ $1,$ ..., $k.
$ If we identify the vertex $v_{1}$ and $v_{2},$ then we have $(K_{1}$ $/$ $%
\{v_{1},$ $v_{2}\})$ $\#_{v_{1}\#v_{2}}$ ... $\#_{v_{1}\#v_{2}}$ $(K_{k}$ $/$
$\{v_{1},$ $v_{2}\})$ of type $L_{k}.$ So, if we construct the edge-iterated
glued graph of a finite directed graph $G,$ with its gluing components which
are generated by edges like above $K_{j}$'s, then we can prove the above
theorem. Again, notice that if $G$ is graph-isomorphic to the edge-iterated
glued graph $K$ with its gluing components $K_{1},$ ..., $K_{N},$ then the
total boundary quotient graph $G$ $/$ $\partial $ satisfies that

\strut 

\begin{center}
$G$ $/$ $\partial =\left( K_{1}\text{ }/\text{ }\partial _{1}\right)
\#_{v_{\partial }}$ ... $\#_{v_{\partial }}$ $\left( K_{N}\text{ }/\text{ }%
\partial _{N}\right) ,$
\end{center}

\strut 

where $\partial _{1},$ ..., $\partial _{N}$ are total boundaries of $K_{1},$
..., $K_{N}$, respectively, and where $v_{\partial }$ is the base point $%
v_{\partial _{1}}$ $\#$ ... $\#$ $v_{\partial _{N}}.$ 

\strut

\strut \strut

\strut

\subsection{Admissible Boundary Quotient Graphs}

\strut

\strut

In this section, we will define and observe the admissible boundaries of
finite simplicial directed graphs and the corresponding boundary quotient
graphs.

\strut \strut

\begin{definition}
Let $G$ be a finite simplicial directed graph and let $\partial _{a}$ be a
boundary of $G.$ The boundary $\partial _{a}$ is said to be an admissible
boundary if $\partial _{a}$ is the boundary of $G$ satisfying the following
conditions;

\strut 

(1)\ there is no admissible finite paths connecting $v_{1}$ and $v_{2},$ for
all pair $(v_{1},$ $v_{2})$ $\in $ $\partial _{a}$ $\times $ $\partial _{a}$
such that $v_{1}$ $\neq $ $v_{2}.$

\strut 

(2) the set $\partial _{a}$ is the maximal subset of $V(G)$ satisfying the
condition (1).
\end{definition}

\strut

We can get that the boundary quotient $\partial _{a}$ is an invariant on
finite simplicial directed graphs with its mixed maximal types. First let's
show the following lemma.

\strut

\begin{lemma}
Let $G$ be a finite simplicial connected directed graph which is
graph-isomorphic to its maximal $CT$-iterated glued graph $K$ of its mixed
maximal type $((C_{N_{1}},$ ..., $C_{N_{s}}),$ $(T_{1},$ ..., $T_{r}),$ $(1,$
..., $1)).$ If $G_{1},$ ..., $G_{r}$ are full subgraphs of $G$ with their
type $T_{1},$ ..., $T_{r},$ respectively, as the gluing components of $K,$
then the admissible boundary $\partial _{a}$ of $G$ is contained in $\cup
_{j=1}^{r}$ $V(G_{j}).$ In particular, $\partial _{a}$ $=$ $\cup _{j=1}^{r}$ 
$(V(G_{j})$ $\cap $ $\partial _{a}).$
\end{lemma}

\strut

\begin{proof}
Remark that the given graph $G$ is connected. Since $G$ is finite
simplicial, there is a unique the maximal $CT$-iterated glued graph $K$
graph-isomorphic to $G.$ Let $K_{1},$ ..., $K_{n}$ be the gluing components
of $K$. Then each $K_{j}$ is graph-isomorphic to a full subgraph of $G$ and
it is of type $C_{N_{j}}$ or of type $T_{j}$ or of type $1.$ Assume that $%
K_{j}$ is of type $C_{N_{j}}.$ Then, for any pair $(v_{1},$ $v_{2})$ in $%
V(K_{j})$ $\times $ $V(K_{j}),$ there always exists a finite path connecting 
$v_{1}$ and $v_{2}$ in $K_{j},$ because $K_{j}$ is a circulant graph. So, $%
\partial _{a}$ $\cap $ $V(K_{j})$ $=$ $\varnothing ,$ for such $j.$ Now,
suppose that $K_{j}$ is of type $T.$ Then $\partial _{a}$ $\cap $ $K_{j}$ is
either empty or non-empty (See the following example). Therefore,

\strut

$\ \ \ \ \ \ \ \ \ \ \ \ \ \ \ \ \ \ \ \ \ \ \ \ \ \partial _{a}$ $=$ $\cup
_{i=1}^{r}\left( \partial _{a}\text{ }\cap \text{ }V(K_{j_{i}})\right) ,$

\strut

where $K_{j_{i}}$ is of type $T_{j_{i}},$ for all $i$ $=$ $1,$ ..., $r.$
\end{proof}

\strut

The above lemma shows how we can determine the admissible boundaries for the
simplicial connected directed graphs.\strut  The admissible boundary $%
\partial _{a}$ of a finite simplicial connected directed graph $G$ is the
subset of the disjoint union of the vertex sets of gluing components of type 
$T.$

\strut

\begin{example}
Consider the following three non-isomorphic trees $T_{1},$ $T_{2}$ and $T_{3}
$, where

\strut 

$\ \ \ \ V(T_{1})=\{v_{1}^{1},v_{2}^{1},v_{3}^{1}\}$ and $%
E(T_{1})=\{e_{1}^{1}=v_{1}^{1}e_{1}^{1}v_{2}^{1},$ $%
e_{2}^{1}=v_{3}^{1}e_{2}^{1}v_{1}^{1}\}$

$\ \ \ \ V(T_{2})=\{v_{1}^{2},v_{2}^{2},v_{3}^{2}\}$ and $%
E(T_{2})=\{e_{1}^{2}=v_{1}^{2}e_{1}^{2}v_{2}^{2},$ $%
e_{2}^{2}=v_{1}^{2}e_{2}^{2}v_{3}^{2}\}$

$\ \ \ \ V(T_{3})=\{v_{1}^{3},v_{2}^{3},v_{3}^{3}\}$ and $%
E(T_{3})=\{e_{1}^{3}=v_{2}^{3}e_{1}^{3}v_{1}^{3},$ $%
e_{2}^{3}=v_{3}^{3}e_{3}^{2}v_{1}^{3}\}.$

\strut 

Then the admissible boundaries $\partial _{a1},$ $\partial _{a2}$ and $%
\partial _{a3}$ of $G_{1},$ $G_{2}$ and $G_{3}$ are

\strut 

$\ \ \ \ \ \ \ \ \ \ \ \ \partial _{a1}=\varnothing ,$ \ \ $\partial
_{a2}=\{v_{2},$ $v_{3}\}$ \ \ and \ \ $\partial _{a3}=\{v_{2},$ $v_{3}\}.$

\strut 

So, the admissible boundary quotient graph $T_{1}$ $/$ $\partial _{a1}$ $=$ $%
T_{1}$ and the admissible boundary quotient graphs $T_{2}$ $/$ $\partial
_{a2}$ and $T_{3}$ $/$ $\partial _{a3}$ are graphs with

\strut \strut 

$\ \ \ \ \ \ \ \ \ V\left( T_{2}\text{ }/\text{ }\partial _{a2}\right)
=\{v_{1}^{2},$ $v_{\partial _{a2}}\}$ and $E\left( T_{2}\text{ }/\text{ }%
\partial _{a2}\right) =\{e_{1},e_{2}\},$

\strut 

where $e_{1}=v_{1}^{2}e_{1}v_{\partial _{a2}}$ and $e_{2}=v_{1}^{2}e_{2}v_{%
\partial _{a2}},$ and

\strut 

$\ \ \ \ \ \ \ \ \ V\left( T_{3}\text{ }/\text{ }\partial _{a3}\right)
=\{v_{1}^{3},v_{\partial _{a3}}\}$ and $E\left( T_{3}\text{ }/\text{ }%
\partial _{a3}\right) =\{f_{1},$ $f_{2}\},$

\strut 

where $f_{1}=v_{\partial _{a3}}$ $f_{1}$ $v_{1}^{3}$ and $f_{2}$ $=$ $%
v_{\partial _{a_{3}}}$ $f_{2}$ $v_{1}^{3}.$
\end{example}

\strut

In the above example, we can observe that the admissible boundary quotient
graphs $T_{2}$ $/$ $\partial _{a2}$ and $T_{3}$ $/$ $\partial _{a3}$ are
graph-isomorphic, via the graph isomorphism $g$ $:$ $T_{2}$ $/$ $\partial
_{a2}$ $\rightarrow $ $T_{3}$ $/$ $\partial _{a3}$ mapping $v_{1}^{2}$ $%
\mapsto $ $v_{\partial _{a3}}$ and $v_{\partial _{a2}}$ $\mapsto $ $%
v_{1}^{3}.$ So, unfortunately, the admissible boundary quotient is not an
invariants on the finite simplicial directed graphs.

\strut

\begin{remark}
The admissible boundary quotient is not an invariant on finite directed
trees and hence it is not an invariant on finite simplicial connected
directed graphs.
\end{remark}

\strut \strut \strut

Recall that we say a finite simplicial directed graph $G$ is of maximal type
if it is graph-isomorphic to the maximal $CT$-iterated glued graph $K$ and
each gluing component is of type $C_{N}$ or of type $T$ or of type $1.$ The
given graph $G$ is said to be of \textbf{mixed maximal type}, if there
exists at least one distinct pair $(K_{1},$ $K_{2})$ of the gluing
components such that $K_{1}$ and $K_{2}$ have different types. i.e., if $%
K_{1}$ is of type $C_{N}$ (or $T$), then $K_{2}$ is of type $T^{\prime }$
(resp. $C_{M}$). As we have seen before, the admissible boundary quotient $%
\partial _{a}$ on finite simplicial directed graph is not an invariant.
However, we can get the following result;

\strut

\begin{theorem}
The admissible boundary quotient $\partial _{a}$ is an invariant on finite
simplicial \textbf{connected} directed graphs of \textbf{mixed maximal type}%
. i.e., if $G_{1}$ and $G_{2}$ are graph-isomorphic finite simplicial
directed graphs and if $\partial _{a1}$ and $\partial _{a2}$ are
corresponding admissible boundaries of $G_{1}$ and $G_{2},$ respectively,
then the boundary quotient graphs $G_{1}$ $/$ $\partial _{a1}$ and $G_{2}$ $/
$ $\partial _{a2}$ are also graph-isomorphic. And the converse is also true.
\end{theorem}

\strut

\begin{proof}
($\Rightarrow $) Suppose $G_{1}$ and $G_{2}$ are graph-isomorphic and assume
that $g$ $:$ $G_{1}$ $\rightarrow $ $G_{2}$ is the graph-isomorphism. Then,
by definition, the map $g$ is a bijection between $V(G_{1})$ and $V(G_{2}),$
preserving the admissibility on $G_{1}$. Take the admissible boundaries $%
\partial _{a1}$ and $\partial _{a2}$ of $G_{1}$ and $G_{2},$ respectively.
Since $g$ preserves the admissibility, we can get that $g\left( \partial
_{a1}\right) $ $=$ $\partial _{a2},$ by the maximality of $\partial _{a1}$
and $\partial _{a2}.$ Therefore, $G_{2}$ $/$ $\partial _{a2}$ $=$ $g\left(
G_{1}\right) $ $/$ $g\left( \partial _{a1}\right) .$ They are
graph-isomorphic via $g^{\symbol{94}}$ $:$ $G_{1}$ $/$ $\partial _{a1}$ $%
\rightarrow $ $G_{2}$ $/$ $\partial _{a2},$ where

\strut

\ $\ \ \ \ \ \ \ \ \ \ \ \ \ \ \ \ \ g^{\symbol{94}}(v)$ $=$ $\left\{ 
\begin{array}{lll}
g(v) &  & \text{if }v\in V(G_{1})\text{ }\setminus \text{ }\partial _{a1} \\ 
v_{\partial _{a2}} &  & \text{if }v=v_{\partial _{a1}}.
\end{array}
\right. $

\strut

Remark that $V\left( G_{1}\text{ }/\text{ }\partial _{a1}\right) =\left(
V(G_{1})\setminus \partial _{a1}\right) \cup \{v_{\partial _{a1}}\}.$ Since $%
g(\partial _{a1})$ $=$ $\partial _{a2},$ the isomorphism $g^{\symbol{94}}$
is well-determined by $g.$ (Notice that we use neither the assumption that $%
G_{1}$ and $G_{2}$ are of mixed maximal type nor $G_{1}$ and $G_{2}$ are
connected. See the next proposition.)

\strut

($\Leftarrow $) Assume that $G_{1}$ and $G_{2}$ are not graph-isomorphic.
Recall that every finite simplicial directed graph is graph-isomorphic to
its unique maximal $CT$-iterated glued graph. Take the finite directed
graphs $K_{1}$ and $K_{2}$ which are the maximal $CT$-iterated glued graphs
of $G_{1}$ and $G_{2},$ respectively. Suppose that $K_{j}$ has its gluing
components $K_{j1},$ ..., $K_{jn_{j}},$ for $j$ $=$ $1,$ $2.$ Since $G_{1}$
and $G_{2}$ are not graph-isomorphic, $K_{1}$ and $K_{2}$ are also not
graph-isomorphic. In other words, $K_{1}$ and $K_{2}$ have different types.
Assume that the gluing components $K_{1r_{i}}$ are of type $T_{r_{i}}$ in $%
\{K_{11},$ ..., $K_{1n_{1}}\}$ and $K_{2s_{i}}$ are of type $T_{s_{i}}$ in $%
\{K_{21},$ ..., $K_{2n_{2}}\}.$ By the previous lemma, we have that

\strut

$\ \ \ \ \ \ \ \ \ \ \ \ \ \ \ \ \ \ \ \partial _{aj}=$ $\ \underset{i}{\cup 
}\left( \partial _{a}\cap V(K_{jr_{i}})\right) ,$ for $j$ $=$ $1,$ $2.$

\strut

In other words, the admissible boundary quotient $\partial _{aj}$ does not
act on the gluing components of type $C_{N}$. So, the admissible boundary
quotient graphs of $G_{1}$ and $G_{2}$ are not graph-isomorphic, since $K_{1}
$ and $K_{2}$ have the different type with same number of vertices.
\end{proof}

\strut

In the previous theorem, we show that if two finite simplicial connected
directed graphs with the mixed maximal types are graph-isomorphic, then
their admissible boundary quotient graphs are also graph-isomorphic. In the
next proposition, we will consider the general case when we just have two
graph-isomorphic finite simplicial directed graphs (which are not
necessarily of mixed maximal type). We can easily verify that their
admissible boundary quotient graphs are also graph-isomorphic, by ($%
\Rightarrow $) in the proof of the previous theorem. In ($\Rightarrow $) of
the previous proof, we did not use the assumption that $G_{1}$ and $G_{2}$
are of mixed maximal type. Therefore, by ($\Rightarrow $) of the previous
proof, we have;

\strut

\begin{proposition}
Let $G_{j}$ be a finite simplicial directed graph (not necessarily connected
or mixed maximal type) and let $\partial _{aj}$ be the admissible boundary
of $G_{j},$ for $j$ $=$ $1,$ $2.$ If $G_{1}$ and $G_{2}$ are
graph-isomorphic, then $G_{1}$ $/$ $\partial _{a1}$ and $G_{2}$ $/$ $%
\partial _{a2}$ are graph-isomorphic. $\square $
\end{proposition}

\strut \strut

Again, remark that the converse of the previous proposition does not hold
true, by the previous example. The previous theorem provides the condition
which makes the converse of the previous proposition hold true. The
condition we found is when $G_{1}$ and $G_{2}$ are connected and they are of
mixed maximal type.

\strut

\begin{example}
(1) Let $G$ $=$ $G_{1}$ $\#_{v_{3}}$ $G_{2}$ be a glued graph of $G_{1}$ and 
$G_{2}$, with its glued vertex $v_{3},$ where $G_{1}$ is the graph of type $%
C_{3}$ with $V(G_{1})$ $=$ $\{v_{1},$ $v_{2},$ $v_{3}\}$ and $E(G_{1})$ $=$ $%
\{e_{1},$ $e_{2},$ $e_{3}\}$, and $G_{2}$ is the graph of type $T,$ with $%
V(G_{2})$ $=$ $\{v_{3},$ $v_{4},$ $v_{5}\}$ and $E(G_{2})$ $=$ $\{e_{4}$ $=$ 
$v_{3}$ $e_{4}$ $v_{4},$ $e_{5}$ $=$ $v_{3}$ $e_{5}$ $e_{5}\}.$ i.e., the
graph $G$ has the maximal $CT$-iterated glued graph of \textbf{mixed}
maximal type $(C_{3},$ $T).$ Then the admissible boundary $\partial _{a}$ $=$
$\{v_{4},$ $v_{5}\}.$ So, we have the admissible quotient graph $G$ $/$ $%
\partial _{a}$ having its glued graph $K_{1}$ $\#_{v_{3}}$ $K_{2},$ where $%
K_{1}$ is the full subgraph $G_{1}$ in $G$ and $K_{2}$ is the graph with $%
V(K_{2})$ $=$ $\{v_{\partial _{a}},$ $v_{3}\}$ and $E(K_{2})$ $=$ $\{f_{1},$ 
$f_{2}\},$ where $f_{1}$ $=$ $v_{3}$ $f_{1}$ $v_{\partial _{a}}$ and $f_{2}$ 
$=$ $v_{3}$ $f_{2}$ $v_{\partial _{a}}.$

\strut 

(2) Now, let $G^{\prime }$ $=$ $G_{1}^{\prime }$ $\#_{v_{3}}$ $G_{2}^{\prime
}$ be a glued graph of $G_{1}^{\prime }$ and $G_{2}^{\prime }$, with its
glued vertex $v_{3},$ where $G_{1}^{\prime }$ is the graph of type $C_{3}$
with $V(G_{1}^{\prime })$ $=$ $\{v_{1},$ $v_{2},$ $v_{3}\}$ and $%
E(G_{2}^{\prime })$ $=$ $\{e_{1},$ $e_{2},$ $e_{3}\}$, and $G_{2}^{\prime }$
is the graph of type $T^{\prime },$ with $V(G_{2}^{\prime })$ $=$ $\{v_{3},$ 
$v_{4},$ $v_{5}\}$ and $E(G_{2}^{\prime })$ $=$ $\{e_{4}$ $=$ $v_{4}$ $e_{4}$
$v_{3},$ $e_{5}$ $=$ $v_{5}$ $e_{5}$ $e_{3}\}.$ i.e., the graph $G^{\prime }$
has the maximal $CT$-iterated glued graph of \textbf{mixed} maximal type $%
(C_{3},$ $T^{\prime }).$ Then the admissible boundary $\partial _{a}^{\prime
}$ $=$ $\{v_{4},$ $v_{5}\}.$ So, we have the admissible quotient graph $G$ $/
$ $\partial _{a}^{\prime }$ having its glued graph $K_{1}^{\prime }$ $%
\#_{v_{3}}$ $K_{2}^{\prime },$ where $K_{1}^{\prime }$ is the full subgraph $%
G_{1}^{\prime }$ in $G^{\prime }$ and $K_{2}$ is the graph with $V(K_{2})$ $=
$ $\{v_{\partial _{a}^{\prime }},$ $v_{3}\}$ and $E(K_{2})$ $=$ $\{f_{1},$ $%
f_{2}\},$ where $f_{1}$ $=$ $v_{\partial _{a}^{\prime }}$ $f_{1}$ $v_{3}$
and $f_{2}$ $=$ $v_{\partial _{a}^{\prime }}$ $f_{2}$ $v_{3}.$

\strut 

(3) The gluing components $G_{2}$ of $G$ in (1) and $G_{2}^{\prime }$ of $%
G^{\prime }$ in (2) are graph-isomorphic, by the previous example. However,
the admissible boundary quotient graphs $G$ $/$ $\partial _{a}$ in (1) and $%
G^{\prime }$ $/$ $\partial _{a}^{\prime }$ $\ $in (2) are not
graph-isomorphic.

\strut 

(4) It is easy to check that the one-flow circulant graphs $C_{N_{1}}$ and $%
C_{N_{2}}$ are graph-isomorphic if and only if $N_{1}$ $=$ $N_{2}$ in $\Bbb{N%
}$ $\setminus $ $\{1\}.$ Moreover, the admissible boundaries $\partial _{1}$
and $\partial _{2}$ of $C_{N_{1}}$ and $C_{N_{2}}$ are empty. Therefore the
admissible boundary quotient graphs of them are $C_{N_{1}}$ and $C_{N_{2}},$
respectively. Thus, anyway, the boundary quotient is an invariant on finite
one-flow circulant graphs.
\end{example}

\strut \strut

\strut

\strut

\section{Subgraph Boundary Quotinet Graphs and Subgraph Boundary Index}

\strut

\strut

Let $G$ be a finite directed graph such that $\left| V(G)\right| $ $<$ $%
\infty $ and $\left| E(G)\right| $ $<$ $\infty $ (not necessarily simplicial
and connected). In this chapter, we will define the graph index $Ind_{G}(H)$
of the graph $G$ with respect to its full subgraph $H.$ To do that we will
define the subgraph boundary quotient graph of $G$ $/$ $\partial _{H},$
where $\partial _{H}$ is the subgraph boundary of $H.$

\strut

\begin{definition}
Let $G$ be an arbitrary finite directed graph and let $H$ be a full subgraph
of $G.$ Define the subgraph boundary $\partial _{H}$ by the set $V(H)$ $\cup 
$ $E(H).$
\end{definition}

\strut

Now, we will define the subgraph boundary quotient graph $G$ $/$ $\partial
_{H}.$

\strut

\begin{definition}
Let $G$ be a finite directed graph and $H,$ a full subgraph and let $%
\partial _{H}$ be the subgraph boundary. Define the subgraph boundary
quotient graph $G$ $/$ $\partial _{H},$ by the directed graph with

\strut 

$\ \ \ \ \ \ \ \ \ \ \ \ \ \ \ \ \ \ \ V(G$ $/$ $\partial _{H})$ $=$ $%
\{v_{\partial _{H}}\}$ $\cup $ $\left( V(G)\text{ }\setminus \text{ }%
V(H)\right) .$

and

$\ \ \ \ \ \ \ \ \ \ \ \ \ \ \ \ \ \ \ \ \ \ \ \ \ \ \ E(G$ $/$ $\partial
_{H})$ $=$ $E(G)$ $\setminus $ $E(H),$

\strut 

with the subgraph boundary quotient; if $x_{1}$ $\neq $ $x_{2}$ in $\partial
_{H},$ then identify $v_{\partial _{H}}$ $=$ $x_{1}$ $=$ $x_{2},$ for all
such pair $(x_{1},$ $x_{2})$ in $\partial _{H}$ $\times $ $\partial _{H}.$
\end{definition}

\strut \strut \strut \strut 

\begin{definition}
By $\left| K\right| ,$ we will denote the size $\left| V(K)\text{ }\cup 
\text{ }E(K)\right| $ of the set $V(K)$ $\cup $ $E(K)$ for all finite
directed graphs $K.$ Let $G$ be a finite directed graph and $H,$ a full
subgraph and let $\partial _{H}$ be the subgraph boundary and $G$ $/$ $%
\partial _{H}$, the corresponding subgraph boundary quotient graph of $G$
with respect to $H.$ The number $\exp $ $(\left| G\text{ }/\text{ }\partial
_{H}\right| $ $-$ $1)$ $\equiv $ $e^{\left| G\text{ }/\text{ }\partial
_{H}\right| -1}$ is called the subgraph boundary index of $G$ with respect
to $H,$ and it is denoted by $ind_{G}(H).$
\end{definition}

\strut

Let $K$ be a finite directed graph. Then, since $K$ is finite the cardinality

$\strut $

\begin{center}
$\left| K\right| $ $=$ $\left| V(K)\cup E(K)\right| =\left| V(K)\right|
+\left| E(K)\right| <\infty .$
\end{center}

\strut \strut

If $G$ is a finite directed graph and $H$ is a full subgraph of $G$ and if $G
$ $/$ $\partial _{H}$ is the corresponding subgraph boundary quotient graph,
then

\strut

(2.1) $\ \ \ \ \ \ \ \ \ \ \ Ind_{G}(H)\overset{def}{=}\exp \left( \left| G%
\text{ }/\text{ }\partial _{H}\right| -1\right) \leq \exp \left| G\right|
<\infty .$

\strut

More generally, we can get the following proposition;

\strut

\begin{proposition}
Let $G$ be a finite directed graph and $H,$ a full subgraph of $G$. Then the
subgraph boundary index $Ind_{G}(H)$ $=$ $\frac{Ind_{G}(1)}{Ind_{H}(1)},$
where $1$ means the trivial graph (i.e., $1$ is the one-vertex-no-edge
graph).
\end{proposition}

\strut

\begin{proof}
By definition, $Ind_{G}(H)=\exp \left( \left| G\text{ }/\text{ }\partial
_{H}\right| -1\right) ,$ where $\partial _{H}$ is the subgraph boundary. It
suffices to show that $\left| G\text{ }/\text{ }\partial _{H}\right| $ $=$ $%
\left| G\right| $ $-$ $\left| H\right| +1.$ It is trivial by definition of
the subgraph admissible quotient graph $G$ $/$ $\partial _{H}.$ Observe that

\strut

$\ \ \ 
\begin{array}{ll}
\left| G\text{ }/\text{ }\partial _{H}\right|  & 
\begin{array}{l}
=\left| V(G\text{ }/\text{ }\partial _{H})\right| +\left| E(G\text{ }/\text{ 
}\partial _{H})\right|  \\ 
\end{array}
\\ 
& 
\begin{array}{l}
=\left( \left| V(G)\right| -\left| V(H)\right| +\left| \{v_{\partial
_{H}}\}\right| \right) +\left( \left| E(G)\right| -\left| E(H)\right|
\right)  \\ 
\end{array}
\\ 
& 
\begin{array}{l}
=\left( \left| V(G)\right| +\left| E(G)\right| \right) -\left( \left|
V(H)\right| +\left| E(H)\right| \right) +1 \\ 
\end{array}
\\ 
& =\left| G\right| -\left| H\right| +1.
\end{array}
$

\strut

Thus we have that $\left| G\text{ }/\text{ }\partial _{H}\right| -1=\left|
G\right| -\left| H\right| .$ By taking exponential on both sides, we have

\strut

$\ \ \ \ \ \ \ \ \ \ \ \ \ \ \ \ \ \ Ind_{G}(H)=\exp \left( \left| G\right|
-\left| H\right| \right) =\frac{\exp \left| G\right| }{\exp \left| H\right| }%
.$

\strut

For the trivial full subgraph $1,$ By the very definition, $G$ $/$ $\partial
_{1}$ $=$ $G$ and hence

\strut

$\ \ \ \ \ \ \ \ \ \ \ \ \ \ \ Ind_{G}(1)=\exp \left( \left| G\text{ }/\text{
}\partial _{1}\right| -1\right) =e^{-1}\exp \left| G\right| ,$

\strut

where $\partial _{1}$ is the subgraph boundary of $1.$ Similarly, $%
Ind_{H}(1) $ $=$ $e^{-1}$ $\exp $ $\left| H\right| .$ Therefore, we can get
that

\strut

$\ \ \ \ \ \ \ \ \ \ \ \ \ \ \ \ \ \ \ \ \ \ \ \ \ \ \ \ \ \ \ \ Ind_{G}(H)=%
\frac{Ind_{G}(1)}{Ind_{H}(1)}.$

\strut
\end{proof}

\strut

By the previous proposition, we can get the following simple results;

\strut

\begin{corollary}
Let $G$ be a finite directed graph and $H,$ a full subgraph in $G.$

\strut 

\strut (1) $Ind_{G}(H)=1$ if and only if $H=G.$

(2) $Ind_{G}(H)>1$ if and only if $H$ is properly contained in $G$.
\end{corollary}

\strut

\begin{proof}
The subgraph boundary index $Ind_{G}(H)$ is the number $\exp $ $(\left| G%
\text{ }/\text{ }\partial _{H}\right| $ $-$ $1),$ where $\partial _{H}$ is
the subgraph boundary of $G$ with respect to $H.$ By the previous
proposition, $Ind_{G}(H)$ $=$ $\frac{Ind_{G}(1)}{Ind_{H}(1)},$ where $1$
means the trivial graph.

\strut

(1) Suppose that $Ind_{G}(H)=1.$ Equivalently,

\strut

$\ \ \ \ \ \ \ \ \ \ \ \ \ \ \ \ \ Ind_{G}(1)$ $=$ $Ind_{H}(1)%
\Longleftrightarrow \exp \left| G\right| =\exp \left| H\right| .$

\strut

So, we have that $\left| G\right| =\left| H\right| .$ Since $H$ is a full
subgraph of $G,$ $V(H)$ $\subseteq $ $V(G)$ and $E(H)$ $\subseteq $ $E(G).$
Thus the condition says that

\strut

$\ \ \ \ \ \ \ \ \ \ \ \ \ \ \ \ \ \ \ \ \left| V(G)\right| +\left|
E(G)\right| =\left| V(H)\right| +\left| E(H)\right| .$

\strut

By the disjointness of the vertex set and the edge set, $\left| V(H)\right| $
$=$ $\left| V(G)\right| $ and $\left| E(H)\right| $ $=$ $\left| E(G)\right|
, $ and hence $V(H)$ $=$ $V(G)$ and $E(H)$ $=$ $E(G).$ Therefore, $H$ $=$ $%
G. $ Conversely, if $H$ $=$ $G,$ then

\strut

$\ \ \ \ \ \ \ \ \ \ \ \ \ \ \ \ \ \ \ \ \ \ \ \ \ Ind_{G}(H)$ $=$ $%
Ind_{G}(G)$ $=$ $\frac{Ind_{G}(1)}{Ind_{G}(1)}$ $=$ $1.$

\strut

(2) Suppose that $Ind_{G}(H)$ $>$ $1.$ Then, by the previous proposition, we
have that

\strut

\ \ \ \ \ \ \ $Ind_{G}(H)=\frac{Ind_{G}(1)}{Ind_{H}(1)}>1\Leftrightarrow
Ind_{G}(1)>Ind_{H}(1)\Leftrightarrow \left| G\right| >\left| H\right| .$

\strut

Thus $H$ is proper full subgraph of $G.$ The converse clearly holds true.
\end{proof}

\strut

\begin{proposition}
Let $G$ be a finite directed graph and $H_{1}$ and $H_{2},$ full subgraphs.
Then $\left| H_{1}\right| $ $=$ $\left| H_{2}\right| $ if and only if $%
Ind_{G}(H_{1})$ $=$ $Ind_{G}(H_{2}).$
\end{proposition}

\strut

\begin{proof}
($\Rightarrow $) Let $G$, $H_{1}$ and $H_{2}$ be given as above. By the
previous proposition,

\strut

(2.2)$\ \ \ \ \ \ \ \ \ Ind_{G}(H_{1})=\frac{Ind_{G}(1)}{Ind_{H_{1}}(1)}=%
\frac{Ind_{G}(1)}{Ind_{H_{2}}(1)}=Ind_{G}(H_{2}).$

\strut

We can get that $Ind_{H_{j}}(1)=\exp \left( \left| H_{j}\right| -1\right) ,$
for $j$ $=$ $1,$ $2.$ Since $\left| H_{1}\right| =\left| H_{2}\right| ,$ the
subgraph boundary indices $Ind_{H_{1}}(1)$ and $Ind_{H_{2}}(1)$ coincide.
Thus the second equality of (2.2) holds true.

\strut

($\Leftarrow $) Trivial, by the definition of the subgraph boundary quotient
index.
\end{proof}

\strut

Let $G_{1}$ and $G_{2}$ be finite directed graphs and assume that they are
graph-isomorphic via the graph-isomorphism $g$ $:$ $G_{1}$ $\rightarrow $ $%
G_{2}.$ Let $H_{1}$ be a full subgraph of $G_{1}.$ Then the image $g(H_{1})$
of $H_{1}$ is also a full subgraph of $G_{2}.$ The next theorem shows that
the subgraph boundary index $Ind$ is preserved up to graph-isomorphisms.

\strut

\begin{theorem}
Let $G$ be a finite directed graph and $H,$ a full subgraph. Suppose that
the graph $G^{\prime }$ is graph-isomorphic to $G.$ Then $Ind_{G}(H)$ $=$ $%
Ind_{G^{\prime }}(H^{\prime }),$ where $H^{\prime }$ is the image of $H,$ of
the corresponding graph-isomorphism, in $G^{\prime }$
\end{theorem}

\strut

\begin{proof}
Let $g$ $:$ $G$ $\rightarrow $ $G^{\prime }$ be a graph-isomorphism. Then $g$
preserves the vertex set and the admissibility of $G.$ So, if $H$ is a full
subgraph of $G,$ then the image $H^{\prime }$ $=$ $g(H)$ of $G^{\prime }$ is
also a full subgraph and moreover $H^{\prime }$ is graph-isomorphic to $H,$
via $g^{-1}\mid _{H^{\prime }}.$ Therefore, $\left| H\right| $ $=$ $\left|
H^{\prime }\right| .$ By the little modification of the previous
proposition, $Ind_{G}(H)$ $=$ $Ind_{G^{\prime }}(H^{\prime }).$
\end{proof}

\strut

By the previous theorem, if $G_{1}$ and $G_{2}$ are finite graphs and if
there exists a graph-homomorphism $g$ $:$ $G_{1}$ $\rightarrow $ $G_{2}$
such that (i) $g\left( V(G_{1})\right) $ $\subseteq $ $V(G_{2})$ and (ii) $g$
preserves the admissibility of $G_{1}$ in $G_{2}.$ Then we can regard the
image $g(G_{1})$ of $G_{2}$ as the full subgraph of $G_{2}.$ So, we can
define the index of $Ind_{G_{2}}(G_{1})$ by the subgraph boundary index $%
Ind_{G_{2}}\left( g(G_{1})\right) .$ i.e.,

\strut

\begin{center}
$Ind_{G_{2}}(G_{1})\overset{def}{=}\exp \left( \left| G_{2}\text{ }/\text{ }%
\partial _{g(G_{1})}\right| -1\right) .$
\end{center}

\strut

\begin{definition}
Let $G_{1}$ and $G_{2}$ be finite directed graphs. Then the boundary index $%
Ind_{G_{2}}(G_{1})$ is defined by

\strut 

$\ \ \ Ind_{G_{2}}(G_{1})\overset{def}{=}\left\{ 
\begin{array}{ll}
Ind_{G_{2}}\left( g(G_{1})\right)  & \text{if }\exists \text{ homomorphism }%
g:G_{1}\rightarrow G_{2} \\ 
0 & \text{otherwise.}
\end{array}
\right. $
\end{definition}

\strut

By definition, we have the following theorem.

\strut

\begin{theorem}
Let $G_{1}$ and $G_{2}$ be finite directed graphs. Then $Ind_{G_{2}}(G_{1})$ 
$=$ $1$ $=$ $Ind_{G_{1}}(G_{2})$ if and only if $G_{1}$ and $G_{2}$ are
graph-isomorphic.
\end{theorem}

\strut

\begin{proof}
($\Rightarrow $) Assume that $Ind_{G_{2}}(G_{1})=1.$ This means that there
exists a graph-homomorphism $g$ $:$ $G_{1}$ $\rightarrow $ $G_{2}$ such that
the boundary index $Ind_{G_{2}}(G_{1})$ is the subgraph boundary index $%
Ind_{G_{2}}$ $\left( g(G_{1})\right) $, defined by $\exp $ $((G_{2}$ $/$ $%
\partial _{g(G_{1})})$ $-$ $1)$ $=$ $1.$ Since $Ind_{G_{2}}$ $\left(
g(G_{1})\right) $ $=$ $\frac{Ind_{G_{2}}(1)}{Ind_{g(G_{1})}(1)},$ if this
quantity is $1,$ then $Ind_{G_{2}}(1)$ $=$ $Ind_{g(G_{1})}(1),$ and hence $%
G_{2}$ $=$ $g(G_{1}).$ Equivalently, the graph-homomorphism $g$ is a
graph-isomorphism from $G_{1}$ to $G_{2}.$ Also, by defining the
graph-isomorphism $g^{-1}$ $:$ $G_{2}$ $\rightarrow $ $G_{1},$ we can get
that $Ind_{G_{1}}$ $\left( g^{-1}(G_{2})\right) $ $=$ $1.$ So, $%
Ind_{G_{1}}(G_{2})$ $=$ $1,$ too.

\strut

($\Leftarrow $) Assume that $G_{1}$ and $G_{2}$ are graph-isomorphic with
its graph-isomorphism $g$ $:$ $G_{1}$ $\rightarrow $ $G_{2}.$ Then the image 
$g(G_{1})$ is the full subgraph of $G_{2},$ moreover $g(G_{1})$ $=$ $G_{2}.$
Therefore,

\strut

$\ \ \ \ \ \ \ \ \ \ \ \ \ \ \ \ \ \ \ \ \ \ \ \frac{Ind_{G_{2}}(1)}{%
Ind_{g(G_{1})}(1)}$ $=\frac{\exp \left( \left| G_{2}\right| -1\right) }{\exp
\left( \left| g(G_{1})\right| -1\right) }=1,$

\strut

and hence $Ind_{G_{2}}\left( g(G_{1})\right) =1.$ So, $Ind_{G_{2}}(G_{1})$ $=
$ $1.$ Also, $Ind_{G_{1}}(G_{2})$ $=$ $1,$ via the graph-isomorphism $g^{-1}.
$
\end{proof}

\strut

We can verify the range of the subgraph boundary index $Ind_{G}(.)$ of the
graph $G,$ as a fundtion defined on the set of all full subgraphs of $G.$

\strut

\begin{proposition}
The image of the subgraph boundary index $Ind_{G}(\cdot ),$ as a
full-subgraph-valued function, is contained in the closed interval $[1,$ $%
e^{\left| G\right| -1}].$ $\square $
\end{proposition}

\strut

More explicitly, we have the following theorem.

\strut

\begin{theorem}
The image of the subgraph boundary index $Ind_{G}(\cdot )$ is contained in $%
\{e^{\left| G\right| -k}$ $:$ $k$ $=$ $1,$ ..., $\left| G\right| \}.$
\end{theorem}

\strut

\begin{proof}
By definition, if $H$ is a full subgraph of $G,$ then $Ind_{G}(H)$ $=$ $\exp 
$ $(\left| G\text{ }/\text{ }\partial _{H}\right| $ $-$ $1).$ Futhermore, $%
Ind_{G}(H)$ $=$ $\frac{Ind_{G}(1)}{Ind_{H}(1)}$ $=$ $\frac{e^{\left|
G\right| -1}}{e^{\left| H\right| -1}}.$ Consider a full subgraph $L$ such
that $V(L)$ $=$ $\{v\}$ and $E(L)$ $=$ $\{l\},$ where $l$ is the loop-edge
concentrated on $v.$ Then $Ind_{L}(1)$ $=$ $e^{\left| L\right| -1}$ $=$ $e.$
Now, let $K$ be a full subgraph with $V(K)$ $=$ $\{v_{1},$ $v_{2}\}$ and $%
E(K)$ $=$ $\{e\},$ where $e$ is the edge connection $v_{1}$ and $v_{2},$
with direction. Then $Ind_{K}(1)$ $=$ $e^{\left| K\right| -1}$ $=$ $e^{2}.$
If $H$ is $G,$ itself, then we have $Ind_{G}(H)$ $=$ $1.$ Also, if $H$ is
trivial, then $Ind_{G}(H)$ $=$ $e^{\left| G\right| -1}.$
\end{proof}

\strut \strut

Let's consider the chain of full subgraphs in the given graph $G.$ We denote
the relation [$H$ is a full subgraph of $G$] by [$H$ $<$ $G$]. The finite
inclusions

\strut

\begin{center}
$K_{1}<K_{2}<...<K_{n}<G$
\end{center}

\strut

is called the chain of full subgraphs.

\strut

\begin{proposition}
Let $G$ be a finite directed graph and let $K$ $<$ $H$ $<$ $G$ be a chain of
full subgraphs of $G.$ Then $Ind_{G}(K)$ $=$ $Ind_{G}(H)$ $\cdot $ $%
Ind_{H}(K).$
\end{proposition}

\strut

\begin{proof}
We have that

\strut

$\ \ \ \ \ \ \ \ \ 
\begin{array}{ll}
Ind_{G}(K) & =\frac{Ind_{G}(1)}{Ind_{K}(1)}=\left( \frac{Ind_{G}(1)}{%
Ind_{H}(1)}\right) /\left( \frac{Ind_{H}(1)}{Ind_{K}(1)}\right) \\ 
&  \\ 
& =Ind_{G}(H)\cdot Ind_{H}(K).
\end{array}
$

\strut
\end{proof}

\strut

More generally, we have the following corollary.

\strut

\begin{corollary}
Let $G$ be a finite directed graph and let $K_{1}$ $<$ ... $<$ $K_{n}$ $<$ $G
$ be a chain of full subgraphs. Then

\strut 

$\ \ \ \ \ \ \ \ \ \ \ \ \ \ \ Ind_{G}(K_{1})$ $=$ $\Pi _{j=2}^{n+1}$ $%
Ind_{K_{j}}(K_{j-1}),$

\strut 

with $K_{n+1}$ $\overset{def}{=}$ $G.$ $\square $
\end{corollary}

\strut

Let $K_{1}$ $<$ $K_{2}$ $<$ ... $<$ $K_{n}$ $<$ $G$ be a chain of full
subgraphs. Then we have the following dual structure of it;

\strut

\begin{center}
$(G$ $/$ $\partial _{K_{n}})$ $<$ $(G$ $/$ $\partial _{K_{n-1}})$ $<$ ... $<$
$(G$ $/$ $\partial _{K_{2}})$ $<$ $(G$ $/$ $\partial _{K_{1}}).$
\end{center}

\strut

This is called the dual chain of the given chain.

\strut

\begin{lemma}
Let $K$ $<$ $H$ $<$ $G$ be a chain of full subgraphs. Then $(G$ $/$ $%
\partial _{H})$ $<$ $(G$ $/$ $\partial _{K}).$
\end{lemma}

\strut

\begin{proof}
Clearly, $Ind_{G}(H)$ $\leq $ $Ind_{G}(K).$ Equivalently, $Ind_{H}(1)$ $\geq 
$ $Ind_{K}(1).$ Also since $\partial _{K}$ $\subset $ $\partial _{H},$ the
subgraph boundary quotient graphs $G$ $/$ $\partial _{H}$ and $G$ $/$ $%
\partial _{K}$ satisfy the full-subgraph-inclusion $(G$ $/$ $\partial _{H})$ 
$<$ $(G$ $/$ $\partial _{K}).$
\end{proof}

\strut

The above lemma shows that the dual chain of the chain of full subgraphs is
well-defined, as a chain of full subgraphs.

\strut

\begin{proposition}
Let $K<H<G$ be a chain of full subgraphs. Then

\strut \strut

$\ \ \ \ \ \ \ \ \ \ \ \ \ \ \ \ \ \ \ \ \ Ind_{G\text{ }/\text{ }\partial
_{K}}$ $\left( G\text{ }/\text{ }\partial _{H}\right) $ $=Ind_{H}(K).$
\end{proposition}

\strut

\begin{proof}
By the previous lemma, we have the full-subgraph-inclusion $(G$ $/$ $%
\partial _{H})$ $<$ $(G$ $/$ $\partial _{K}),$ as the dual chain of $K$ $<$ $%
H$ $<$ $G.$ The subgraph boundary index $Ind_{G\text{ }/\text{ }\partial
_{K}}$ $\left( G\text{ }/\text{ }\partial _{H}\right) $ is determined by

\strut

$\ \ Ind_{G\text{ }/\text{ }\partial _{K}}\left( G\text{ }/\text{ }\partial
_{H}\right) =\frac{Ind_{G\text{ }/\text{ }\partial _{K}}(1)}{Ind_{G\text{ }/%
\text{ }\partial _{H}}(1)}$ $=$ $\frac{\exp (\left| G\text{ }/\text{ }%
\partial _{K}\right| -1)}{\exp (\left| G\text{ }/\text{ }\partial
_{H}\right| -1)}$

\strut

$\ \ \ \ \ \ \ \ \ \ \ \ \ \ \ =\frac{Ind_{G}(K)}{Ind_{G}(H)}=\left( \frac{%
Ind_{G}(1)}{Ind_{K}(1)}\right) $ $/$ $\left( \frac{Ind_{G}(1)}{Ind_{H}(1)}%
\right) =\left( \frac{Ind_{G}(1)}{Ind_{K}(1)}\right) \left( \frac{Ind_{H}(1)%
}{Ind_{G}(1)}\right) $

\strut

$\ \ \ \ \ \ \ \ \ \ \ \ \ \ \ =\frac{Ind_{H}(1)}{Ind_{K}(1)}=Ind_{H}(K).$

\strut \strut

Therefore, $Ind_{G\text{ }/\text{ }\partial _{K}}\left( G\text{ }/\text{ }%
\partial _{H}\right) =Ind_{H}(K).$
\end{proof}

\strut

By the previous proposition, generally, we can get that;

\strut

\begin{corollary}
Let $K_{1}$ $<$ ... $<$ $K_{n}$ $<$ $G$ be a chain of full subgraphs and let 
$(G$ $/$ $\partial _{K_{n}})$ $<$ ... $<$ $(G$ $/$ $\partial _{K_{1}})$ be
the corresponding dual chain. Then

\strut

$\ \ \ \ \ \ \ \ \ \ \ \ \ \ \ \ \ \ \ Ind_{G\text{ }/\text{ }\partial
_{K_{i}}}\left( G\text{ }/\text{ }\partial _{K_{j}}\right) $ $=$ $%
Ind_{K_{j}}(K_{i}),$

\strut

for all $i$ $\leq $ $j$ in $\{1,$ ..., $n\}.$ $\square $
\end{corollary}

\strut

\begin{remark}
As we observed before, we can apply the above results to the general chain
of graphs. Assume that we have a chain of graphs $G_{1}$ $<$ $...$ $<$ $%
G_{n+1}.$ Here $G_{i}$ $<$ $G_{i+1}$ means that there exists a
graph-homomorphism $g_{i}$ $:$ $G_{i}$ $\rightarrow $ $G_{i+1}$ such that $%
g_{i}(G_{i})$ is a full subgraph of $G_{i+1},$ for all $i$ $=$ $1,$ ..., $n.$
Then we have its dual chain $(G_{n+1}$ $/$ $\partial _{g_{n}(G_{n})})$ $<$
... $<$ $(G_{n+1}$ $/$ $\partial _{g_{1}(G_{1})})$ $<$ $(G_{n+1}$ $/$ $%
\partial _{i(G_{1})}),$ where $i$ $:$ $G_{1}$ $\rightarrow $ $G_{1}$ is the
identity graph-isomorphism. Then we have the following boundary index
relations;

\strut

(1) \ $Ind_{G_{n+1}}(G_{1})$ $=$ $\Pi _{j=2}^{n+1}$ $Ind_{G_{j}}(G_{j-1}).$

\strut

(2) \ \ $Ind_{G_{n+1}\text{ }/\text{ }\partial _{G_{i}}}\left( G_{n+1}\text{ 
}/\text{ }\partial _{G_{j}}\right) $ $=$ $Ind_{G_{j}}(G_{i}),$ for all $i$ $%
\leq $ $j$ in $\{1,$ ..., $n$ $+$ $1\}.$
\end{remark}

\strut

\strut \strut

\strut \strut

\textbf{References}

\strut

\strut

\begin{quotation}
\strut {\small [1] \ A. G. Myasnikov and V. Shpilrain (Editors), Group
Theory, Statistics, and Cryptography, Contemp. Math 360, AMS (2004)}

{\small [2] \ \ B. Bollobas and E. Szemeredi, Girth of Sparse Graphs, J.
Graph Theory 39 (2002), 194 - 200.}

{\small [3] \ \ D. V. Osin, Algebraic Entropy of Elementary Amenable Groups,
(2005), Preprint.}

{\small [4] \ \ F. Buckley and M. Lewinter, A Friendly Introduction to Graph
Theory, Pearson Education Inc. (2003) ISBN 0-13-066949-0.}

{\small [5] \ \ I. Cho, Random Variables in a Graph }$W^{*}${\small %
-Probability Space, Ph. D. Thesis, Univ. of Iowa (2005)}

{\small [6] \ \ I. Cho, Graph Measure Theory (2005), Preprint.}

{\small [7] \ \ I. Cho, Free Probability on Graph Measure Algebras, in
Progress.}

{\small [8] \ \ J. Lauri and R. Scapellato, Topics in Graph Automorphisms
and Reconstruction, London.Math.Soc.Stu.Text 54, Cambridge Univ. Press (2003)%
}

{\small [9] \ \ M. P. Bellon and C. M. Viallet, Algebraic Entropy, (1998),
Preprint.}

{\small [10] R. Gilman, V. Shpilrain and A. G. Myasnikov (Editors),
Computational and Statistical Group Theory, Contemp. Math. 298, AMS (2002)}

{\small [11] T. W. Hungerford, Algebra (1980)}
\end{quotation}

\end{document}